\documentclass[12pt,psfig]{article}
\usepackage{amsmath,amssymb,epsfig, epstopdf}
\newtheorem{Theorem}{Theorem}[section]
\newtheorem{theorem}{Theorem}[section]
\newtheorem{definition}[theorem]{Definition}

\newtheorem{lemma}[theorem]{Lemma}
\newtheorem{corollary}[theorem]{Corollary}

\def\E{\mathbb{E}}
\def\P{\mathbb{P}}

\def\P{\mathbb{P}}
\def\E{\mathbb{E}}

\def \sS {{\cal S}}

\makeatletter 
\@addtoreset{equation}{section}

\makeatother 
\begin{document}
\title{A Multi Period Equilibrium Pricing Model\thanks{  Work supported by NSERC grants 371653-09, 88051 and MITACS grants 5-26761, 30354 and the Natural Science Foundation of China (10901086).
{\bf{Acknowledgement}}: We would like to thank an anonymous referee and Tom Hurd for helpful comments.}}

\author{\normalsize    Traian A.~Pirvu \\[8pt]
        \small Dept of Mathematics \& Statistics\\
        \small McMaster University \\
        \small 1280 Main Street West \\
        \small Hamilton, ON, L8S 4K1\\
        \small tpirvu@math.mcmaster.ca
        \and
        \normalsize Huayue Zhang \\[8pt]
        \small Dept of Finance\\
        \small  Nankai University\\
        \small 94 Weijin Road \\
        \small Tianjin, China, 300071 \\
        \small  hyzhang69@nankai.edu.cn
\vspace*{0.8cm}}

\maketitle

\baselineskip=12pt
\begin{abstract}

In this paper, we propose an equilibrium pricing model in a dynamic
multi-period stochastic framework with uncertain income streams. 
In an incomplete market, there exist two traded
risky assets (e.g. stock/commodity and weather derivative) and a
non-traded underlying (e.g. temperature).  The risk preferences are
of exponential (CARA) type with a stochastic coefficient of risk
aversion. Both time consistent and time inconsistent
trading strategies are considered. We obtain the equilibriums
prices of a contingent claim written on the risky asset and
non-traded underlying.  By running numerical experiments we
examine how the equilibriums prices vary in response to changes in model
parameters.

\end{abstract}

\textbf{Keywords:}\ Time inconsistent control, incomplete market,
equilibrium price.
\newpage{}

\section{Introduction}

Hitherto, there has been an increasing literature on pricing
contingent claims written on non-tradable underlyings in a dynamic
multi-period equilibrium framework. One example of such contingent
claim is a weather derivative, in which case the underlying is the
temperature process. One approach in pricing this financial
instruments is to use a multi-period stochastic equilibrium model.
In financial economics there is a big literature on this issue.

Rubinstein (1976) considers a multi-period state-preference
equilibrium model without explicit modelling of
production/investment. Brennan (1979) looks at a multiperiod equilibrium problem in which the representative
agent exhibits constant risk aversion. Bhattacharya (1981) extends
the model of Rubinstein (1976) to show that risk/return tradeoffs are
linear relations linking instantaneous expected returns of assets to
the instantaneous covariances of returns with aggregate consumption.
Bizid and Jouini (2001) derives restrictions on the equilibrium
state-price deflators independent on the choices of utility function
in an incomplete market. Camara (2003) obtains preference-free
option prices in a discrete equilibrium model where representative
agent has exponential utility and aggregate wealth together with the
underlying asset price have transformed normal distributions.

Our paper presents a partial equilibrium model with two exogenous
assets, one tradable and one non-tradable. A derivative security,
written on the tradable and non-tradable assets, is priced in equilibrium
by a representative agent who receives unspanned random income within an
incomplete multi-period market.  Cao and Wei (2004), Lee and Oren
(2009), Lee and Oren (2010), Cheridito  et al (2011)  are related to our work. 
Cao and Wei (2004) generalizes the model of Lucas (1978) to
provide an equilibrium framework for valuing weather derivatives in
a multi-period setting. Lee and Oren (2009) explores a single-period
equilibrium pricing model in a multi commodity setting and
mean-variance preferences. Lee and Oren (2010) is a follow up in a
multi-period framework. Cheridito  et al (2011) establishes results
on the existence and uniqueness of equilibrium in dynamically
incomplete financial markets with preferences of monetary type and
heterogeneous agents.

In our model the representative agent has risk preferences of
exponential type, time and state dependent, with time changing coefficient of risk aversion.
Inspired by Gordon and St-Amour (2000), we assume that
the risk aversion coefficient is a stochastic process. Gordon and St-Amour (2000)
motivates this change by the fact that it can explain asset-price movements which fixed preference
paradigm can not explain. Lately, the issue of time changing
risk aversion received some attention in the financial literature.
For instance,  Barberis (2001) considers a model in which the loss aversion depends on prior gains and losses, so it
 may change through time. Danthine et al (2004) allows the representative agent's coefficient of
 relative risk aversion to vary with the underlying economy's growth rate.
Gordon and St-Amour (2004) explains equity premium puzzle by
state-dependent risk preferences. Yuan and Chen (2006) shows
 that dynamic risk aversion plays a critical role in the
 dynamics of asset price fluctuations.

 A time changing risk preference leads to time inconsistent investment
strategies. It means that an investor may have an incentive to deviate
 from the optimal strategies which he/she computed at some past time.
 In order to overcome this issue, Bjork and Murgoci (2010) develops a theory for
 stochastic control problems which are time inconsistent in the sense that they do not admit a
 Bellman optimality principle; they consider the subgame perfect Nash
 equilibrium strategies as a substitute for the optimal strategies.

  This paper considers both time consistent and time inconsistent optimal strategies. 
  Time consistent optimal strategies are the subgame perfect strategies. Time inconsistent
  optimal strategies are the classical optimal strategies given that the agent does not update
  his/her risk preferences. Time consistent (inconsistent) equilibrium price is defined by imposing
  the market clearing condition for time consistent (inconsistent) optimal strategies.

 In the present work the exogenous assets have stochastic drifts and volatilities.
 It may be that they depend on each other (in a weather model we consider a commodity whose
 volatility is influenced by the temperature process). A derivative security is priced in equilibrium
 within this model. Our main result is an iterative algorithm which yields the equilibrium prices (time consistent and time inconsistent). At each stage the equilibrium prices depend on the current risk aversion level, and all previous wealth and risk aversion levels. The algorithm constructs recursively one period pricing
 kernels. Moreover, the time inconsistent equilibrium pricing kernel equals the marginal utility. 
 We prove that the equilibrium pricing measures are martingale measures so the equilibrium prices
 (time consistent and time inconsistent) are arbitrage free. Numerical experiments shed light into the structure of equilibriums prices. We show that utility indifference prices and equilibriums ones
  are different. The utility indifference price was introduced by Hodges and Neuberger (1989). By now, there are several papers on this topic; we recall only a few. \footnote{ on discrete time Musiela and Zariphopoulou (2004), Musiela et al (2009); on continuous time Henderson (2002), Musiela and Zariphopoulou (2004); for an overview see Henderson Hobson (2004).} Pirvu and Zhang (2011) derives utility indifference prices in a model with time changing risk aversion.

  Next, we consider examples in which the non-traded underlying affects the
  income, thus creating an incentive for the agent to hold the derivative in order to hedge the risk. Our plots show that equilibriums prices are increasing in risk aversion, fact explained by an increased hedging demand.
  We add a stochastic unspanned component to the income and this slightly decreases the equilibrium
  prices. This is explained by a decrease in the hedging demand, which is a consequence of income being only partial affected by the non-traded underlying. Finally, in a regime switching model we explore the effect of changing risk aversion. Here we find that an increase of twenty percent in risk aversion can cause a
  percentage change in the time consistent equilibrium price (with the benchmark being the time inconsistent equilibrium price) anywhere between minus four and seventeen percent.

The optimal strategies are obtained by backward induction. First order conditions together with the market clearing give the equilibrium prices. We consider a partial equilibrium model because of the problem we want to address. However our method can be easily extended to a full equilibrium model (in fact the computations are simpler in that case).

The remainder of this paper is organized as follows. Section 2
presents the model. Section 3 provides the equilibrium pricing valuation. Numerical experiments are presented in the section 4.
Proofs of the results are delegated to an appendix.

\section{The Model}
We consider a multi-period stochastic model of investment.  The
trading horizon is $[0,T],$ with $T$ a exogenous finite horizon.
There are $N+1$ trading dates: $t_{n}=nh,$ for $n=0,1,\cdots,N,$
and $h=\frac{T}{N}.$ Let
$(b^1,b^2,...,b^d):=(b_{t_n}^1,b_{t_n}^2,...,b_{t_n}^d)_{n=0,1,...,\infty},$
be a $d$-dimensional binomial random walk on a complete probability
space $(\Omega,\mathcal {F},\{{\mathcal {F}}_{t_n}\}, \mathbb{P})$.
The random walk is assumed symmetric under $\mathbb{P}$ in the sense
that
\begin{equation}\label{*0}
\P(\Delta b_{t_n}^i=\pm 1)=1/2,\qquad i=1,2,...,d.
\end{equation}
 There are three securities available for trading in our model; a riskless bond, a primary asset (e.g. stock
 or commodity) and a derivative security. We take the bond as numeraire, thus it can be assumed to offer zero
interest rate. The primary asset price process $C:=\{C_{t_n} ;n=0,1,\ldots
,N\}$, follows the difference equation :
\begin{eqnarray}\label{c}
\left\{\begin{array}{ll}\Delta C_{t_n}=C_{t_n}
(\mu^c_{t_n}h+\sigma^c_{t_n}\sqrt{h} \Delta
b_{t_n}^1),\ \ n=0,1,...N-1\\\
\\C_0=c>0,\label{st1}\qquad\qquad\end{array} \right.
\end{eqnarray}
for some adapted drift process $\mu^c:=\{\mu_t^c;t=0,h\ldots,
(N-1)h,Nh\}$ and volatility process
$\sigma^c:=\{\sigma_t^c;t=0,h\ldots, (N-1)h,Nh\}$ which are chosen
so that the commodity price remains positive. The derivative
security is written on the primary asset and a non-tradable underlying. 
The non-traded asset $S:=\{S_{t_n} ;n=0,1,\ldots
,N\}$, follows the difference equation :
\begin{eqnarray}
\left\{\begin{array}{ll} \Delta S_{t_n}=S_{t_n}
(\mu^s_{t_n}h+\sigma^s_{t_n}\sqrt{h} (\rho\Delta
b_{t_n}^1+\sqrt{1-\rho^2}\Delta b_{t_n}^2)),\ \ n=0,1,...,N-1
\\S_0=s>0,\label{st2}\qquad\qquad\end{array}
\right.
\end{eqnarray}for some adapted drift process $\mu^s:=\{\mu_t^s;t=0,h\ldots,
(N-1)h,Nh\},$ volatility process
$\sigma^s:=\{\sigma_t^s;t=0,h\ldots, (N-1)h,Nh\},$ and a correlation
coefficient $\rho,$ with $|\rho|<1.$ The two dimensional process
$P=(C, S),$ exogenously given, is referred to as the forward
process.

The derivative security $D:=\{D_{t_n} ;n=0,1,\ldots
,N\}$ is to be priced in equilibrium. Therefore we have a partial
equilibrium model. It is motivated by a situation in which the
primary asset and derivative are priced in different markets. Take
for example energy and weather derivatives. Although correlated (in
California it was observed a high correlation between energy prices
and temperature process) energy and weather derivatives are priced within different markets. 
\subsection{Trading strategies}
 Let $\alpha_{t_n}$ be the wealth invested in the primary asset
at time $t_n,$ and  $\beta_{t_n}$ the number of shares of
derivative held at time  $t_n$; denote $\pi_{t_n}:=\{\alpha_{t_n},\beta_{t_n}\}\in
{\mathcal{F}_{t_n}}, n=0,1\ldots, N-1.$  The value of a
self-financed portfolio satisfies the following stochastic
difference equation:
\begin{eqnarray}
\Delta
X^{\pi}_{t_n}&:=&\alpha_{t_n}(\mu^c_{t_n}h+\sigma^c_{t_n}\sqrt{h}
\Delta b_{t_n}^1)+\beta_{t_n} (\Delta
D_{t_{n}}+\varphi(t_{n},C_{t_{n}},S_{t_{n}})h). \label{dxt}
\end{eqnarray}
Here $\varphi(t_{n},C_{t_{n}},S_{t_{n}})h$ is the dividend paid for holding the stock on $[t_n , t_{n+1}].$ At maturity, $t_N:=T=Nh,$ the representative agent in this economy receives
 random income $I_{t_N},$ which is ${\mathcal{F}_T}$
adapted. Thus, his/her final wealth is
\begin{equation}\label{wt}
W_{t_N}^\pi=X_{t_N}^{\pi}+I_{t_N}.
\end{equation}

The random income may depend on all the random walks
$\{b^1,b^2,...,b^d\},$ so it may not be spanned by the existing
assets. This makes our market model incomplete.

\subsection{Risk Preferences}

The representative agent utility is assumed to be of exponential
type, time and state dependent. Moreover, the coefficient of absolute
risk aversion is a stochastic process $\gamma_{t_n},$ $\{{\mathcal {F}}_{t_n}\}$
adapted,  $n=0,1\ldots, N-1.$ More precisely

\begin{equation}\label{riskp}
U(x,t_n,\omega)=-\exp(-\gamma_{t_n} (\omega)x).
\end{equation}
This modeling approach is not
new; in the introduction we pointed out a number of papers which
consider time changing risk aversion. The performance of an
investment strategy $\pi$ is measured by the above expected utility
criterion applied to the final wealth. At time $t_n$ the optimization criterion is 

\begin{equation}\label{val}
 \sup_{\pi \in \Pi_{t_{n}} } \E[- e^{-\gamma_{t_n}
 W^{\pi}_{t_N} } | {\mathcal {F}}_{t_n}  ].
 \end{equation}
Here $n=0,1\ldots, N-1,$ $W^{\pi}_{t_N}$ is given by \eqref{wt}, and $\Pi_{t_{n}}$ denotes the
set of admissible trading strategies,
\begin{equation}\label{dd}
\Pi_{t_{n}}:=\{\pi_{t_n}, \pi_{t_{n+1}}, \cdots,
\pi_{t_{N-1}}:\pi_{t_k}\in{\mathcal{F}}_{t_{k}}, \ {\rm{such\
that}}\,\,\E |X^{\pi}_{t_k}|<\infty, k=n,n+1\ldots N-1 \}.
\end{equation}

\subsection{Optimal Time Inconsistent Strategies}
They are the classical optimal strategies given that the risk preferences are not
updated. More precisely $\hat{\pi}\in \Pi_{t_{0}} $ is an optimal time inconsistent strategy
if it satisfies 

\begin{eqnarray}\label{*0} { \hat{\pi} }=\arg\sup_{\pi\in\Pi_{t_0}}\E^{\P}
[-\exp(-\gamma_{t_0} W^{\pi}_{t_{N}})|  {\mathcal {F}}_{t_0}   ].\label{v}\end{eqnarray}
They are called time inconsistent because they fail to remain
optimal at later times $t_n,$ in the sense that

\begin{eqnarray}\label{*00}  { \hat{\pi} }\neq\arg\sup_{\pi\in\Pi_{t_n}}\E^{\P}
[-\exp(-\gamma_{t_n}  W^{\pi}_{t_{N}})|  {\mathcal {F}}_{t_n} ].\label{v}\end{eqnarray}

\subsection{Optimal Time Consistent Strategies}

In this section, we introduce the optimal time consistent strategies defined as subgame perfect strategies.
First, let us consider the time period $[(N-1)h,Nh]$ (recall that
$t_{N}=T=Nh$). At time $(N-1)h$ consider the optimization problem:
\begin{eqnarray}  (P1) \qquad 
\sup_{\pi\in\Pi_{t_{N-1}}}
\E [-\exp(-\gamma_{t_{N-1}}  W^{\pi}_{t_N} ) |  {\mathcal {F}}_{t_{N-1}}].
\label{v1}\end{eqnarray} In
our model $\sup$ in (P1) is attained and we denote 

\begin{eqnarray} {\pi}_{t_{N-1}}^*=
\arg\max_{\pi\in \Pi_{t_{N-1}}} \E [-\exp(-\gamma_{t_{N-1}}  W^{\pi}_{t_N} ) |  {\mathcal {F}}_{t_{N-1}}] .
\label{v^11}\end{eqnarray} 
 On the time period $[(N-2)h, Nh]$,
restrict to the trading strategies $\pi$ be of the form:
\begin{eqnarray}
\pi=\left\{\begin{array}{ll}\pi^*_{t_{N-1}},\ \ {\rm{on}}\ [(N-1)h,Nh),\\
\\ \pi_{t_{N-2}},\ \ {\rm{on}}\ [(N-2)h,(N-1)h),\label{alpha2}\qquad\qquad\end{array} \right.
\end{eqnarray}
for an arbitrary ${\mathcal {F}}_{t_{N-2}}-$ adapted control
$\pi_{t_{N-2}}$ such that $(\pi_{t_{N-2}},\pi^*_{t_{N-1}}) \in
\Pi_{t_{N-2}};$ consider the optimization problem
\begin{eqnarray} (P2) \qquad \sup_{\pi \in \Pi_{t_{N-2}}}
\E [-\exp(-\gamma_{t_{N-2}} W^{\pi}_{t_N}) |   {\mathcal {F}}_{t_{N-2}} ].
\label{v123}\end{eqnarray}
In our model $\sup$ in (P2) is attained and we denote 
\begin{eqnarray} (\pi^*_{t_{N-2}}, \pi^*_{t_{N-1}})=
\arg \max_{\pi\in \Pi_{t_{N-2}}}
  \E [-\exp(-\gamma_{t_{N-2}} W^{\pi}_{t_N}) |   {\mathcal {F}}_{t_{N-2}} ].\label{v*11}\end{eqnarray}
 Further we proceed iteratively. On the time period $[(N-n)h,Nh]$ one restricts
to trading strategies $\pi$ of the form:
\begin{equation}
\pi =\begin{cases} \pi^*_{t_k},
\ \ {\rm{for}}\ k=N-(n-1),N-(n-2),\cdots, N-1,\\
\\\pi_{t_k}\ \ {\rm{for}} \,\,\,k=N-n,
\label{alphan}\qquad\qquad\end{cases}
\end{equation}
for an arbitrary ${\mathcal {F}}_{t_{N-n}}-$ adapted control
$\pi_{t_{N-n}}$ such that $(\pi_{t_{N-n}},\pi_{t_k}^*)_{\{k=
N-n+1,\cdots, N-1\} }\in \Pi_{t_{N-n}}.$ Consider the optimization problem 
\begin{eqnarray} (Pn)\qquad 
\max_{\pi \in \Pi_{t_{N-n}}}
\E [-\exp(-\gamma_{t_{N-n}} W^{\pi}_{t_N}) |   {\mathcal {F}}_{t_{N-n}} ]. \label{v123}\end{eqnarray}
The $\sup$ in (Pn) is attained and we denote 
\begin{eqnarray} (\pi^*_{t_{N-n}},
\pi^*_{t_{N-n+1}}, \cdots, \pi^*_{t_{N-1}})=\arg \max_{\pi\in
\Pi_{t_{N-n}}}
 \E [-\exp(-\gamma_{t_{N-n}} W^{\pi}_{t_N}) |   {\mathcal {F}}_{t_{N-n}} ].\label{v*11}\end{eqnarray}
The optimal time consistent strategy is ${\pi^*}= (\pi^*_{t_{0}},
\pi^*_{t_{1}}, \cdots, \pi^*_{t_{N}}).$

\subsection{Time Consistent versus Time Inconsistent Strategies}

Let us recall that time inconsistencies in this model are due to 
time changing risk aversion. Indeed in the case of constant risk
aversion, the optimal time consistent
and time inconsistent strategies coincide, i.e. $\pi^*=\hat{\pi}.$
In general, it is hard to show that the optimal time consistent strategy outperforms the
optimal time inconsistent strategy. We show that this is the case if the risk preferences are updated in the following two
period model. Let us assume that $I_{t_2}=0,$ and only one asset is available
for trading (the primary asset with constant drift and volatility). It is claimed that
\begin{equation}\label{&&}
 \E [-\exp(-\gamma_{t_{1}}W^{\pi^*}_{t_{2}}) | \mathcal{F}_{t_{1}}]\geq  \E [-\exp(-\gamma_{t_{1}}W^{\hat{\pi}}_{t_{2}}) | \mathcal{F}_{t_{1}}].
\end{equation} 
Indeed in this model it can be shown that $\pi^*_{t_{1}}=\hat{\pi}_{t_{1}},$ and hence $W^{\pi^*}_{t_{1}}=W^{\hat{\pi}}_{t_{1}}.$ By the definition of optimal time consistent strategies, \eqref{&&} yields.

\section{Equilibrium Valuation}

We assume that there exists a representative agent with risk preferences given
by \eqref{riskp}. The representative agent trades in $C$ and $D$ in order to maximize
expected utility of his/her final wealth. This can be achieved by the optimal time inconsistent
strategy $\hat{\pi}=(\hat{\alpha}, \hat{\beta})$ if the representative agent does not update his/her
risk preferences. Otherwise, the optimal time consistent strategy  ${\pi}^{*}=({\alpha}^{*}, {\beta}^{*})$ will be used.
Thus, depending on weather or not the representative agent updates his/her risk preferences we introduce two notions
of equilibriums: time consistent and time inconsistent. They are given by the market clearing condition in the formal definition
below.

\begin{definition} Given the terminal payoff $D_{t_N}$ and the dividend stream
$\varphi(t_{n},C_{t_{n}},S_{t_{n}})h,$ $D_{t_n}$ is a time consistent equilibrium price if and only if 
$$\beta_{t_n}^*=1,$$
 for every $n=0,1,...,N-1.$ Likewise,  $D_{t_n}$ is a time inconsistent equilibrium price if and only if 
$$\hat{\beta}_{t_n}=1,$$
 for every $n=0,1,...,N-1.$
\end{definition}
This simply says that there is one unit of derivative in the market
and it is priced such that ``it is optimal'' (time consistent or time inconsistent) for the representative
agent to acquire it. The interest of the representative agent in
holding the derivative comes from the fact that his/her income is
exposed to the non-traded asset $S;$ thus, the risk of income
fluctuations due to $S$ can be hedged by trading $D.$ Let $r^c_{t_{n}},$ given by
\begin{equation}\label{1*}
r^c_{t_{n}}:=\frac{\mu^c_{t_{n}}}{\sigma^c_{t_{n}}},
\end{equation}
be the market price of risk (MPR) for the primary asset which is assumed positive.
 We choose time length $h$ small enough such that $1\geq r_{t_n}^c\sqrt{h}.$

\subsection{Single Period}
In this case the optimal time consistent and time inconsistent strategies coincide and so do the corresponding equilibriums. Define $A_{t_{N-1}}:=\{ \Delta
b_{t_{N-1}}^1=1\}$ and $A_{t_{N-1}}^c:=\{\Delta b_{t_{N-1}}^1=-1\}.$

\begin{Theorem} The equilibrium price
at time $t_{N-1}$ is given by
$$D_{t_{N-1}}=\E ^{Q^*}[D_{t_N}+\varphi(t_{N-1},C_{t_{N-1}},S_{t_{N-1}})h |   {\mathcal {F}}_{t_{N-1}}],$$
where the probability measure $Q^*$ is defined by
$$\frac{dQ^*}{d\P}=\Lambda_{t_N}\E [\frac{dQ^*}{d\P} |   {\mathcal {F}}_{t_{N-1}} ].$$
 The pricing kernel $\Lambda_{t_N}$ is
 \begin{eqnarray}\label{Di} \Lambda_{t_N}:=
\left\{\begin{array}{ll}  \lambda_{t_{N-1}}\frac{e^{-\gamma_{t_{N-1}}
[D_{t_N}+I_{t_N}] }}{\E[e^{-\gamma_{t_{N-1}}
[D_{t_N}+I_{t_N}] }|  A_{t_{N-1}}\vee  {\mathcal {F}}_{t_{N-1}}  ]},\ \ \ {\rm{if}} \
\ \ \omega\in A_{t_{N-1}}
 \\  \lambda_{t_{N-1}}\frac{e^{-\gamma_{t_{N-1}}
[D_{t_N}+I_{t_N}] }}{\E[e^{-\gamma_{t_{N-1}}
[D_{t_N}+I_{t_N}] }|  A^c_{t_{N-1}}\vee  {\mathcal {F}}_{t_{N-1}} ]}  ,\ \ \ {\rm{if}}
\ \ \  \omega\in A_{t_{N-1}}^c,\label{st2}\qquad\qquad\end{array}
\right.
\end{eqnarray}
 
  with
\begin{eqnarray}\label{D}\lambda_{t_{N-1}}=
\left\{\begin{array}{ll} {1-r^c_{t_{N-1}}\sqrt{h}},\ \ \ {\rm{if}} \
\ \ \omega\in A_{t_{N-1}}
 \\{1+r^c_{t_{N-1}}\sqrt{h}},\ \ \ {\rm{if}}
\ \ \  \omega\in A_{t_{N-1}}^c. \label{st2}\qquad\qquad\end{array}
\right.
\end{eqnarray}
The optimal trading strategy is given by

\begin{eqnarray}\label{eqq}
{\alpha}_{t_{N-1}}^*&=&\frac{1}{2\gamma_{t_{N-1}} \sigma^c_{t_{N-1}}\sqrt{h}}
\log\bigg(\frac{1+r^c_{t_{N-1}}\sqrt{h}}{1-r^c_{t_{N-1}}\sqrt{h}}\bigg)\nonumber\\
&+& \frac{1}{2\gamma(J_{t_{N-1}})\sigma^c_{t_{N-1}}\sqrt{h}}\log \bigg(
\frac{\E[  e^{-\gamma_{t_{N-1}}
[D_{t_N}+I_{t_N}] } | A_{t_{N-1}} \vee  {\mathcal {F}}_{t_{N-1}}]}{\E[  e^{-\gamma_{t_{N-1}}
[D_{t_N}+I_{t_N}] }  |
A_{t_{N-1}}^c \vee  {\mathcal {F}}_{t_{N-1}}]}\bigg).\nonumber
\end{eqnarray}

\end{Theorem}

Proof of this Theorem is done in Appendix A.
\begin{flushright}
$\square$
\end{flushright}

Next we prove that the probability measure $Q^*$ is a martingale measure so the equilibrium
prices are arbitrage free.

\begin{lemma}\label{m1}
 If the dividend $\varphi=0,$ then the traded assets $\{ C_{t_{n}}\}_{n=N-1,N}$ and $\{ D_{t_{n}}\}_{n=N-1,N}$ are martingales under $Q^*.$
\end{lemma}

Proof: $\{ D_{t_{n}}\}_{n=N-1,N}$ is martingale under $Q^*$ by definition. Next we show that
$\{ C_{t_{n}}\}_{n=N-1,N}$ is martingale under $Q^*.$
It suffices to prove that with $n=N-1$

\begin{eqnarray}
\E^{Q^{*}}[\frac{C_{t_{n+1}}}{C_{t_n}}|{\mathcal{F}_{t_n}}]&=&
\E^{Q^{*}}[1+\sigma^c_{t_{n}}\sqrt{h}(r_{t_{n}}^c\sqrt{h}+\Delta
b^1_{t_{n}})|{\mathcal{F}_{t_n}}]\nonumber\\
&=&1+\sigma^c_{t_{n}}\sqrt{h}\E^{Q^{*}}[(r_{t_{n}}^c\sqrt{h}+\Delta
b^1_{t_{n}})|{\mathcal{F}_{t_n}}]\\&=&1.
\end{eqnarray}
This is the case if
\begin{equation}\label{oo}\E^{Q^{*}}[(r_{t_{n}}^c\sqrt{h}+\Delta
b^1_{t_{n}})|{\mathcal{F}_{t_n}}]=0.\end{equation}

When $n=N-1$, \eqref{oo} is equivalent to
$$
\frac{r_{t_{N-1}}^c\sqrt{h}+1}{2}\E
\left[\Lambda_{t_{N}}\E [\frac{dQ^*}{d\P}|  {\mathcal {F}}_{t_{N-1}}]
|{A_{t_{N-1}}}\vee  {\mathcal {F}}_{t_{N-1}}\right]+\frac{r_{t_{n}}^c\sqrt{h}-1}{2}\E\left[\Lambda_{t_{N}}\E[\frac{dQ^*}{d\P}|  {\mathcal {F}}_{t_{N-1}}]
|{A^c_{t_{N-1}}}\vee  {\mathcal {F}}_{t_{N-1}}\right]=$$$$
\frac{1-(r^c_{t_{N-1}})^2h}{2}\E [\frac{dQ^*}{d\P}|  {\mathcal {F}}_{t_{N-1}}]+
\frac{(r^c_{t_{N-1}})^2h-1}{2}\E[\frac{dQ^*}{d\P} |  {\mathcal {F}}_{t_{N-1}}]
=0,$$
so the claim yields.

\begin{flushright}
$\square$
\end{flushright}

\subsection{Multiple Periods}
Let us define the sets $A_{t_{N-n}}:=\{ \Delta b_{t_{N-n}}^1=1\}$ and
$A_{t_{N-n}}^c:=\{\Delta b_{t_{N-n}}^1=-1\};$
to ease notations, let us denote
$$\E_{t_n}[\cdot]:=\E[ \cdot|  {\mathcal {F}}_{t_{n}}],\,\, \E_{t_{n}}[\cdot|\mathcal{G}]:= \E[\cdot|\mathcal{G}\vee  {\mathcal {F}}_{t_{n}} ],$$
 for every
$\mathcal{G}\subset\mathcal{F}.$ The equilibrium prices are computed by a recursive
algorithm. Imagine that they were found at all prior times and now we
want to find them at $t_{N-n}.$ Define the random variables
$Y^{*}_{t_{N-(n-1)}}$ and $\hat{Y}_{t_{N-(n-1)}}$ by:
\begin{eqnarray}e^{-\gamma_{t_{N-n}} Y^{*}_{t_{N-(n-1)}}}:&=&
\E_{t_{N-n+1}}\bigg[e^{-\gamma_{t_{N-n}}{[{\displaystyle{\sum_{k=N-(n-1)}^{N-1}}}\Delta
{X}_{t_k}^{*}}+I_{t_N}]} \bigg]\label{ynn},
\end{eqnarray}
and
\begin{eqnarray}e^{-\gamma_{t_{N-n}}\hat{Y}_{t_{N-(n-1)}}}:&=&
\E_{t_{N-n+1}}\bigg[e^{-\gamma_{t_{0}}{[{\displaystyle{\sum_{k=N-(n-1)}^{N-1}}}\Delta
\hat{X}_{t_k}}+I_{t_N}]}\bigg].\label{ynn1}
\end{eqnarray}
Here$$\Delta
{X}_{t_k}^{*}=\alpha_{t_k}^*(\mu_{t_k}^ch+\sqrt{h}\sigma_{t_k}^c\Delta
b_{t_{k}}^1) +\Delta D_{t_k}+\varphi(t_{k},C_{t_{k}},S_{t_{k}})h,$$
for any $k=N-n+1,\cdots, N-2,N-1;$ the optimal time consistent strategy
$\alpha_{t_k}^*$ is given by 
\begin{eqnarray}\label{eu}
\alpha_{t_k}^*&=& \frac{1}{2\gamma_{t_{k}}  \sigma^c_{t_{k}}\sqrt{h}}
\log\bigg(\frac{1+r^c_{t_{k}}\sqrt{h}}{1-r^c_{t_{k}}\sqrt{h}}\bigg)\nonumber\\
&+& \frac{1}{2\gamma_{t_{k}} \sigma^c_{t_{k}}\sqrt{h}}\log \bigg(
\frac{\E_{t_{k}} [  e^{-\gamma_{t_{k}}
[D_{t_{k+1}}+Y^*_{t_{k+1}}] } | A_{t_{k}}]}{\E_{t_{k}}[  e^{-\gamma_{t_{k}}
[D_{t_{k+1}}+Y^*_{t_{k+1}}] }  |
A_{t_{k}}^c]}\bigg)
.\end{eqnarray}
Moreover $$\Delta
\hat{X}_{t_k}=\hat{\alpha}_{t_k}(\mu_{t_k}^ch+\sqrt{h}\sigma_{t_k}^c\Delta
b_{t_{k}}^1) +\Delta D_{t_k}+\varphi(t_{k},C_{t_{k}},S_{t_{k}})h,$$
for any $k=N-n+1,\cdots, N-2,N-1;$ the optimal time inconsistent strategy
$\hat{\alpha}_{t_k}$ is given by by
\begin{eqnarray}\label{eu1}
\hat{\alpha}_{t_k}&=& \frac{1}{2\gamma_{t_{0}}    \sigma^c_{t_{k}}\sqrt{h}}
\log\bigg(\frac{1+r^c_{t_{k}}\sqrt{h}}{1-r^c_{t_{k}}\sqrt{h}}\bigg)\nonumber\\
&+& \frac{1}{2\gamma_{t_{0}}   \sigma^c_{t_{k}}\sqrt{h}}\log \bigg(
\frac{\E_{t_{k}} [  e^{-\gamma_{t_{0}}
[D_{t_{k+1}}+\hat{Y}_{t_{k+1}}] } | A_{t_{k}}]}{\E_{t_{k}}[  e^{-\gamma_{t_{0}}
[D_{t_{k+1}}+\hat{Y}_{t_{k+1}}] }  |
A_{t_{k}}^c]}\bigg).
\end{eqnarray}
Notice that $Y^{*}_{t_{N-(n-1)}}$ and $\hat{Y}_{t_{N-(n-1)}}$ are the certainty equivalents (time consistent and time inconsistent) at
time $n-1.$ In the special case of constant coefficient of absolute risk aversion they are equal. Next, define
the one step period pricing kernels $\Lambda^{*}_{t_{N-n+1}}$
and $\hat{\Lambda}_{t_{N-n+1}}$  by

\begin{eqnarray}\label{lambda} \Lambda^{*}_{t_{N-n+1}}:=
\left\{\begin{array}{ll}  \lambda_{t_{N-n}} \frac{e^{-\gamma_{t_{N-n}}[D_{t_{N-n+1}}+Y^{*}_{t_{N-n+1}}]}}{\E_{t_{N-n}}[ e^{-\gamma_{t_{N-n}}[D_{t_{N-n+1}}+Y^{*}_{t_{N-n+1}}]} | A_{t_{N-n}}   ]}  ,\ \ \ {\rm{if}} \
\ \ \omega\in A_{t_{N-n}}
 \\  \lambda_{t_{N-n}} \frac{e^{-\gamma_{t_{N-n}}[D_{t_{N-n+1}}+Y^{*}_{t_{N-n+1}}]}}{\E_{t_{N-n}}[ e^{-\gamma_{t_{N-n}}[D_{t_{N-n+1}}+Y^{*}_{t_{N-n+1}}]} | A^c_{t_{N-n}}   ]},\ \ \ {\rm{if}}
\ \ \  \omega\in A_{t_{N-n}}^c,\label{st21}\qquad\qquad\end{array}
\right.
\end{eqnarray}

\begin{eqnarray}\label{lambda1} \hat{\Lambda}_{t_{N-n+1}}:=
\left\{\begin{array}{ll}  \lambda_{t_{N-n}} \frac{e^{-\gamma_{t_{0}}[D_{t_{N-n+1}}+\hat{Y}_{t_{N-n+1}}]}}{\E_{t_{N-n}}[ e^{-\gamma_{t_{0}}[D_{t_{N-n+1}}+\hat{Y}_{t_{N-n+1}}]} | A_{t_{N-n}}   ]}  ,\ \ \ {\rm{if}} \
\ \ \omega\in A_{t_{N-n}}
 \\  \lambda_{t_{N-n}} \frac{e^{-\gamma_{t_{0}}[D_{t_{N-n+1}}+\hat{Y}_{t_{N-n+1}}]}}{\E_{t_{N-n}}[ e^{-\gamma_{t_{0}}[D_{t_{N-n+1}}+\hat{Y}_{t_{N-n+1}}]} | A^c_{t_{N-n}}   ]}.\ \ \ {\rm{if}}
\ \ \  \omega\in A_{t_{N-n}}^c,\label{st20}\qquad\qquad\end{array}
\right.
\end{eqnarray}
Here
\begin{eqnarray}\label{D}\lambda_{t_{N-n}}=
\left\{\begin{array}{ll} {1-r^c_{t_{N-n}}\sqrt{h}},\ \ \ {\rm{if}} \
\ \ \omega\in A_{t_{N-n}}
 \\{1+r^c_{t_{N-n}}\sqrt{h}},\ \ \ {\rm{if}}
\ \ \  \omega\in A_{t_{N-n}}^c. \label{st2}\qquad\qquad\end{array}
\right.
\end{eqnarray}

The following Theorem is the main result of the paper.

\begin{Theorem}The time consistent equilibrium price
at time $t_{N-n}$ is given by
$$D_{t_{N-n}}=\E_{t_{N-n}}^{Q^*}[D_{t_{N-n+1}}+\varphi(t_{N-n},C_{t_{N-n}},S_{t_{N-n}})h],$$
 where the probability measure $Q^*$ is defined by
$$\frac{dQ^*}{d\P}=\Lambda^{*}_{t_N}\Lambda^{*}_{t_{N-1}}\ldots\Lambda^{*}_{t_1}.$$
The optimal time consistent strategy (in the primary asset) is ${\alpha^*}= (\alpha^*_{t_{0}},
\alpha^*_{t_{1}}, \cdots, \alpha^*_{t_{N}}),$ with $\alpha_{t_k}^*$ defined by
\eqref{eu} . The time inconsistent equilibrium price at time $t_{N-n}$ is given by
$$D_{t_{N-n}}=\E_{t_{N-n}}^{\hat{Q}}[D_{t_{N-n+1}}+\varphi(t_{N-n},C_{t_{N-n}},S_{t_{N-n}})h],$$
 where the probability measure $\hat{Q}$ is defined by
$$\frac{d\hat{Q}}{d\P}=\hat{\Lambda}^{}_{t_N}\hat{\Lambda}^{}_{t_{N-1}}\ldots\hat{\Lambda}^{}_{t_1}.$$
The optimal time inconsistent strategy (in the primary asset) is ${\hat{\alpha}}= (\hat{\alpha}_{t_{0}},
\hat{\alpha}_{t_{1}}, \cdots, \hat{\alpha}_{t_{N}}),$ with $\hat{\alpha}_{t_k}$ defined by
\eqref{eu1} .  
 \end{Theorem}

Proof of this Theorem is done in Appendix B.

\begin{flushright}
$\square$
\end{flushright}

\noindent For the time inconsistent equilibrium price, we recover the following classical result.

\begin{corollary} The time inconsitent pricing kernel equals the
marginal utility, i.e.,
\begin{eqnarray}
\hat{\Lambda}_{t_{N-n+1}}=\frac{\E_{t_{N-n+1}}[U'(\hat{W}_{t_N})]}{\E_{t_{N-n}}[{U'(\hat{W}_{t_N})}]},
\end{eqnarray}
where  $U(x)=-e^{-\gamma x},$ and  $\hat{W}_{t_N}$ ( see \eqref{wt} ) is the optimal time inconsistent wealth.
\end{corollary}
Proof of this Corollary is done in Appendix C.

\begin{flushright}
$\square$
\end{flushright}

\begin{lemma}\label{m11}
 If the dividend $\varphi=0,$ then the traded assets $\{ C_{t_{n}}\}_{n=0,1,\ldots,N}$ and $\{ D_{t_{n}}\}_{n=0,1,\ldots N}$ are martingales under $Q^*.$
\end{lemma}

Proof: The proof is similar to Lemma \ref{m1} so is skipped.

\begin{flushright}
$\square$
\end{flushright}

\section{Numerical examples}
We specialize to a regime switching model.
A discrete time finite state homogeneous Markov
chain (MC) $J:=(J_{t_n})_{n=0,1,...,\infty}$ is defined on
$(\Omega,\mathcal {F},\{{\mathcal {F}}_{t_n}\}, \mathbb{P})$ and it
takes values in the state space $\sS=\{{\bf 0},{\bf1}\}$ (which represents two states of the market,
bull and bear). The $n-$step transition matrix $P^{(n)}=(P_{ij}^{n}),$ is defined by
$$P^{(n)}_{ij}: =\P(J_{t_n}=j|J_{t_0}=i),\ \ i,j={\bf 0},{\bf1}\qquad n=
0,1,...,\infty,$$ where $P^{(0)}_{ij}=1$ when $i=j,$ otherwise
$P^{(0)}_{ij}=0.$ We assume that the distribution of $J_0$ is known,
and $$\P(J_0=i|{{\mathcal {F}}_{0}})=\P(J_0=i), \qquad i={\bf
0},{\bf1}.$$ The risk aversion is defined by $\gamma_{t_n}= \gamma(J_{t_n}.)$
In this section we give a concrete example; take the electricity
industry one of the most weather-sensitive businesses in the
economy. When the temperature increases there is a higher demand for
electricity due to the usage of air conditioners. In turn this
will lead to higher energy prices. In our model an energy provider
hedges the weather exposure by selling one share of weather
derivative to the representative agent. This is designed such that
it has a higher payoff when temperature is high. The representative
agent has an incentive to buy this product because of his/her income
exposure to weather.

Consider an European call option on the temperature with a strike
price $K=10$; assume that $h=0.3,$ $C_0=c=10,\, S_0=s=10$ (this is normalized and it corresponds to 85 Fahrenheit degrees),$\,
\rho=0.5,\, \mu^c=0.1,\ \sigma^c=0.2,\ \mu^s = 0.3,\,\sigma^s =
0.50.$ The energy price process $(C_t)$ follows (\ref{c}); the
market price of risk (MPR) of the commodity $(r_t^c)$ is defined
$$(r_{t_n}^c)^2=(\arctan(S_t)+\frac{\pi}{2}).$$
Therefore higher temperatures lead to an increased (MPR).

\subsection{Single period}
In this section, we present a numerical example of single-period, $N
=1.$ Recall that in this case time consistent equilibrium coincides with
time inconsistent equilibrium. Let $D_1=(S_1-K)^{+},$ and for simplicity assume $\varphi=0.$

\subsubsection{Equilibrium Price versus Indifference Price}
The paths of the typical trajectories of the forward processes
$(C,S),$ are plotted in Fig.1, and the sample path of MPR in  Fig.2.

It is easy, within our model, to compare numerically equilibrium
price and indifference price of $D_1.$This is done in Fig. 3, Fig. 4
and Fig. 5, where we take $\rho=0.5$ and  $\gamma=0.7.$

Next we introduce the income $I_1 =7e^{-0.5(S_1-s)h}.$ Fig. 6 shows
that the equilibrium price is an increasing function of the risk
aversion $\gamma.$ This is explained by an increase in the hedging
demand when the agent becomes more risk averse. This increase is due
to weather impact on the income.

Next we add a nonspanned component to the income
 $$I_1=7e^{-0.5(S_1-s)h}+5e^{0.1h}1_{\{\Delta b_0^1=1,\Delta
b_0^3=1\}}+4e^{0.1h}1_{\{\Delta b_0^1=1,\Delta
b_0^3=-1\}}+$$$$+2e^{0.1h}1_{\{\Delta b_0^1=-1,\Delta
b_0^3=1\}}+e^{0.1h}1_{\{\Delta b_0^1=-1,\Delta b_0^3=-1\}}.$$\\

\noindent Fig 7. shows the effect of this addition. The equilibrium
price of the derivative is slightly lower with nonspanned income.
This is explained by the fact that in this case only a part of the
income is affected by weather, whence a lower hedging demand.

\subsection{Two periods}
Take $N=2,$ $D_2=(S_2-K)^{+},$ and $\varphi=0.$ Moreover
$$I_{2}= 7e^{-0.5(S_2-s)h}+5e^{0.1h}1_{\{\Delta b_1^1=1,\Delta
b_1^3=1\}}+4e^{0.1h}1_{\{\Delta b_1^1=1,\Delta
b_1^3=-1\}}+$$$$+2e^{0.1h}1_{\{\Delta b_1^1=-1,\Delta
b_1^3=1\}}+e^{0.1h}1_{\{\Delta b_1^1=-1,\Delta b_1^3=-1\}}.$$ Fig. 8
plots the time inconsistent equilibrium price with and without unspanned income. The
effect of the unspanned income becomes more pronounced with $N=2.$

\noindent Fig. 9 plots the effect of time changing risk aversion. We
are interested in the percentage change of the time consistent equilibrium price
when the benchmark is the time inconsistent equilibrium price. We allowed the income to depend on the state of the economy.

$$I_2=7e^{(-0.5(S_1-s))h}+10e^{0.03h}1_{\{J_0=\textbf{0},\Delta
b_1^3=1\}}+8e^{0.03h}1_{\{J_0=\textbf{0},\Delta b_1^3=-1\}}+$$$$
+5e^{0.03h}1_{\{ J_0=\textbf{1},\Delta
b_1^3=1\}}+4e^{0.03h}1_{\{J_0=\textbf{1},\Delta b_1^3=-1\}}.$$

\section{Appendix}
\subsection{Appendix A: Proof of Theorem 3.1}
 From (\ref{dxt}) and (\ref{wt}), it follows that
\begin{eqnarray}\label{4}
\E_{t_{N-1}}[-\exp(-\gamma_{t_{N-1}}(X_{t_N}^{\pi}+I_{t_N}))]
&=&\E_{t_{N-1}}[-\exp(-\gamma_{t_{N-1}}(X_{{t_{N-1}}}+\Delta
X_{t_{N-1}}^{\pi}+I_{t_N}
))]\nonumber\\\notag&=&-e^{-\gamma_{t_{N-1}} x}\E_{t_{N-1}}[\exp(-\gamma_{t_{N-1}}(\Delta
X_{t_{N-1}}^{\pi}+I_{t_N}))]\nonumber\\
&:=&-e^{-\gamma_{t_{N-1}} x}g_{N-1}(\alpha,\beta,\cdot)\nonumber,
\end{eqnarray}
where
\begin{eqnarray}
\Delta
X^{\pi}_{t_{N-1}}:\!\!&\!=\!&\!\!\alpha_{t_{N-1}}(\mu^c_{t_{N-1}}h+\sigma^c_{t_{N-1}}\sqrt{h}
\Delta b_{t_{N-1}}^1)+\beta_{t_{N-1}} (\Delta
D_{t_{N-1}}+\varphi(t_{N-1},C_{t_{N-1}},S_{t_{N-1}})h),\nonumber\\
\label{dxt}
\end{eqnarray}
and
\begin{eqnarray}
&&g_{N-1}(\alpha,\beta,\cdot):=\E_{t_{N-1}}[\exp(-\gamma_{t_{N-1}}(\Delta
X_{t_{N-1}}^{\pi}+I_{t_N}))]\nonumber\\
 &&=\frac{1}{2}e^{-\gamma_{t_{N-1}}\alpha(\mu^c_{t_{N-1}}h+\sigma^c_{t_{N-1}}\sqrt{h}
)}
 \E_{t_{N-1}} \bigg[e^{-\gamma_{t_{N-1}}\beta (\Delta D_{t_{N-1}}
 +\varphi(t_{N-1},C_{t_{N-1}},S_{t_{N-1}})h)}e^{-\gamma_{t_{N-1}} I_{t_N}}| A_{t_{N-1}}\bigg]
\nonumber\\ \ &&+\frac{1}{2}e^{-\gamma_{t_{N-1}}\alpha(\mu^c_{t_{N-1}}h
-\sigma^c_{t_{N-1}}\sqrt{h} )}
 \E_{t_{N-1}} \bigg[e^{-\gamma_{t_{N-1}}\beta(\Delta D_{t_{N-1}}+\varphi(t_{N-1},
 C_{t_{N-1}},S_{t_{N-1}})h)}e^{-\gamma_{t_{N-1}} I_{t_N}}| A_{t_{N-1}}^c\bigg].\nonumber\end{eqnarray}
Recall that $A_{t_{k}}:=\{ \Delta b_{t_{k}}^1=1\}$ and
$A_{t_{k}}^c:=\{\Delta b_{t_{k}}^1=-1\}.$ The function
$g_{N-1}(\alpha,\beta,\cdot) $ has the following properties:
$$g_{N-1}(0,0,\cdot)=\E_{t_{N-1}}[e^{-\gamma_{t_{N-1}}I_{t_N}}]\leq1;$$
For a fixed $\beta,$ it follows that for small $h$
$$g_{N-1}(\infty,\beta,\cdot)=\infty;\  g_{N-1}(-\infty,\beta,\cdot)=\infty;$$
By arbitrage argument it follows that $D_{t_{N-1}}$ belongs to the
interval
$$D_{t_{N-1}}(\omega)\in[\inf_{Q}\E^Q[D_{t_N}],\sup_{Q}\E^Q[D_{t_N}]],$$
where $Q$ ranges over the set of probability measures. Consequently,
$$\Delta D_{t_{N-1}}(\omega)\in[D_{t_N}-\sup_{Q}\E^Q[D_{t_N}],D_{t_N}
-\inf_{Q}\E^Q[D_{t_N}]].$$ Thus, the sets $\{\omega: \Delta D_{t_{N-1}}(\omega)>0\},$ and
$\{\omega :\Delta D_{t_{N-1}}(\omega)<0\}$ have positive probability. This implies that
$$g_{N-1}(\alpha,\infty, \cdot)=\infty;\ g_{N-1}(\alpha,-\infty, \cdot)=\infty.$$
From the above analysis, it follows that the minimum of the function
of $g$ is a critical point. First order conditions lead to
\begin{eqnarray}\label{g1}
\frac{\partial g_{N-1}}{\partial \alpha}
 &=&\frac{-\gamma_{t_{N-1}}(\mu^c_{t_{N-1}}h+\sigma^c_{t_{N-1}}\sqrt{h}
)}{2}\cdot
e^{-\gamma_{t_{N-1}} \alpha(\mu^c_{t_{N-1}}h+\sigma^c_{t_{N-1}}\sqrt{h}
)}\nonumber\\&&\times
 \E_{t_{N-1}} \bigg[e^{-\gamma_{t_{N-1}}\beta (\Delta D_{t_{N-1}}+\varphi(t_{N-1},
 C_{t_{N-1}},S_{t_{N-1}})h)}e^{-\gamma_{t_{N-1}} I_{t_N}}| A_{t_{N-1}}\bigg]
\nonumber\\&+&\frac{-\gamma_{t_{N-1}} (\mu^c_{t_{N-1}}h-\sigma^c_{t_{N-1}}\sqrt{h}
)}{2}\cdot
e^{-\gamma_{t_{N-1}}\alpha(\mu^c_{t_{N-1}}h-\sigma^c_{t_{N-1}}\sqrt{h}
)}\nonumber\\
&&\times
 \E_{t_{N-1}} \bigg[e^{-\gamma_{t_{N-1}}\beta (\Delta D_{t_{N-1}}+
 \varphi(t_{N-1},C_{t_{N-1}},S_{t_{N-1}})h)}e^{-\gamma_{t_{N-1}} I_{t_N}}|
 A_{t_{N-1}}^c\bigg]\nonumber\\
 &=&0,\end{eqnarray}
 and
 \begin{eqnarray}
&&\frac{\partial g_{N-1}}{\partial \beta}
 =\frac{1}{2}
e^{-\gamma_{t_{N-1}}\alpha(\mu^c_{t_{N-1}}h+\sigma^c_{t_{N-1}}\sqrt{h}
)}e^{-\gamma_{t_{N-1}}\beta\varphi(t_{N-1},C_{t_{N-1}},S_{t_{N-1}})h}\nonumber\\&&\times
 \E_{t_{N-1}} \bigg[-\gamma_{t_{N-1}}(\Delta D_{t_{N-1}}+\varphi(t_{N-1},C_{t_{N-1}},S_{t_{N-1}})h)\cdot e^{-\gamma_{t_{N-1}}\beta \Delta D_{t_{N-1}}}
 e^{-\gamma_{t_{N-1}} I_{t_N}}| A_{t_{N-1}}\bigg]
\nonumber\\&&+\frac{1}{2}
e^{-\gamma_{t_{N-1}}\alpha(\mu^c_{t_{N-1}}h-\sigma^c_{t_{N-1}}\sqrt{h}
)}e^{-\gamma_{t_{N-1}}\beta\varphi(t_{N-1},C_{t_{N-1}},S_{t_{N-1}})h}\nonumber\\
&&\times
 \E_{t_{N-1}} \bigg[-\gamma_{t_{N-1}} (\Delta D_{t_{N-1}}+\varphi(t_{N-1},C_{t_{N-1}},
 S_{t_{N-1}})h)\cdot e^{-\gamma_{t_{N-1}}\beta \Delta D_{t_{N-1}}}e^{-\gamma_{t_{N-1}} I_{t_N}}|
 A_{t_{N-1}}^c\bigg]\nonumber\\
 &=&0.\end{eqnarray}
Recall that
\begin{equation}\label{delta1}\Delta
D_{t_{N-1}}+\varphi(t_{N-1},C_{t_{N-1}},S_{t_{N-1}})h=
D_{t_N}-\E^{Q^*}_{t_{N-1}}[D_{t_N}],
\end{equation}
for an equilibrium pricing measure $Q^*$ to be found. Since
$\E^{Q^*}_{t_{N-1}}[D_{t_N}]$ is
${\mathcal{F}_{t_{N-1}}}-$measurable, it follows that
\begin{eqnarray}\label{eqq}
{\alpha}_{t_{N-1}}^*&=&\frac{1}{2\gamma(i)\sigma^c_{t_{N-1}}\sqrt{h}}
\log\bigg(\frac{1+r^c_{t_{N-1}}\sqrt{h}}{1-r^c_{t_{N-1}}\sqrt{h}}\bigg)\nonumber\\
&+& \frac{1}{2\gamma_{t_{N-1}}\sigma^c_{t_{N-1}}\sqrt{h}}\log \bigg(
\frac{\E_{t_{N-1}} [e^{-\gamma_{t_{N-1}}\beta_{t_{N-1}}^*
D_{t_N}}e^{-\gamma_{t_{N-1}} I_{t_N}}| A_{t_{N-1}}]}{\E_{t_{N-1}}
[e^{-\gamma_{t_{N-1}}\beta_{t_{N-1}}^* D_N}e^{-\gamma(i) I_{t_N}}|
A_{t_{N-1}}^c]}\bigg).\nonumber
\end{eqnarray}
By the equilibrium condition $\beta_{t_{N-1}}^*=1$ . This together
with $\frac{\partial g_{N-1}}{\partial \beta}=0$ lead to
$$
 \E_{t_{N-1}} \bigg[(\Delta
D_{t_{N-1}}+\varphi(t_{N-1},C_{t_{N-1}},S_{t_{N-1}})h) e^{-\gamma_{t_{N-1}}
\Delta D_{t_{N-1}}}e^{-\gamma_{t_{N-1}} I_{t_N}}| A_{t_{N-1}}\bigg]=
$$$$-e^{2\gamma_{t_{N-1}} \alpha_{t_{N-1}}^*  \sigma^c_{t_{N-1}}\sqrt{h} }
 \E_{t_{N-1}} \bigg[ (\Delta D_{t_{N-1}}+\varphi(t_{N-1},C_{t_{N-1}},
 S_{t_{N-1}})h) e^{-\gamma_{t_{N-1}}\Delta D_{t_{N-1}}}
 e^{-\gamma_{t_{N-1}} I_{t_N}}|
 A_{t_{N-1}}^c\bigg],$$
and
\begin{equation}\label{326}e^{2\gamma_{t_{N-1}}\alpha_{t_{N-1}}^*\sigma^c_{t_{N-1}}\sqrt{h} }=
\frac{(1+r^c_{t_{N-1}}\sqrt{h})\E_{t_{N-1}} [e^{-\gamma_{t_{N-1}}
D_{t_N}}e^{-\gamma_{t_{N-1}} I_{t_N}}|{A_{t_{N-1}}
}]}{(1-r^c_{t_{N-1}}\sqrt{h})\E_{t_{N-1}} [e^{-\gamma_{t_{N-1}}
D_{t_N}}e^{-\gamma_{t_{N-1}} I_{t_N}}|{ A_{t_{N-1}}^c}]}.\end{equation}
Combing the above equations leads to
\begin{eqnarray}
&&\frac{2}{1-r^c_{t_{N-1}}\sqrt{h}}\E^{Q^*}_{t_{N-1}}[D_{t_N}]
 \nonumber\\&=&
\frac{\E_{t_{N-1}}[D_{t_N} e^{-\gamma_{t_{N-1}} D_{t_N}}e^{-\gamma_{t_{N-1}}
I_{t_N}}|A_{t_{N-1}}]}{\E_{t_{N-1}}[ e^{-\gamma_{t_{N-1}}
D_{t_N}}e^{-\gamma_{t_{N-1}}
I_{t_N}}|A_{t_{N-1}}]}\nonumber\\
&+&\frac{(1+r^c_{t_{N-1}}\sqrt{h})}{(1-r^c_{t_{N-1}}\sqrt{h})}
\frac{\E_{t_{N-1}}[D_{t_N} e^{-\gamma_{t_{N-1}} D_{t_N}}e^{-\gamma_{t_{N-1}}
I_{t_N}}|A_{t_{N-1}}^c]}{\E_{t_{N-1}}[ e^{-\gamma(i)
D_{t_N}}e^{-\gamma_{t_{N-1}} I_{t_N}}|A_{t_{N-1}}^c]}\nonumber.
\end{eqnarray}
This together with (\ref{delta1}) imply that:
\begin{eqnarray}\label{eq}
D_{t_{N-1}}-\varphi(t_{N-1},C_{t_{N-1}},S_{t_{N-1}})h&=&\frac{1-r_{t_{N-1}}^c\sqrt{h}}{2}\frac{\E_{t_{N-1}}
[ D_{t_N}e^{-\gamma_{t_{N-1}} D_{t_N}}e^{\gamma_{t_{N-1}} I_{t_N}}    | A_{t_{N-1}}
]}{\E_{t_{N-1}}
[e^{-\gamma(i) D_{t_N}}e^{\gamma(i) I_{t_N}}  | A_{t_{N-1}} ]}\nonumber\\
&+&\frac{1+r_{t_{N-1}}^c\sqrt{h}}{2} \frac{\E_{t_{N-1}} [
D_{t_N}e^{-\gamma_{t_{N-1}} D_{t_N}}e^{\gamma_{t_{N-1}} I_{t_N}}  | A^c_{t_{N-1}}
]}{\E_{t_{N-1}} [e^{-\gamma_{t_{N-1}} D_{t_N}}e^{\gamma_{t_{N-1}} I_{t_N}}  | A^c_{t_{N-1}} ]}.\nonumber
\end{eqnarray}
Thus, the equilibrium price is
$$D_{t_{N-1}}=\E_{t_{N-1}}[(D_{t_N}+\varphi(t_{N-1},C_{t_{N-1}},S_{t_{N-1}})h)
\Lambda_{t_{N}}],$$
where $\Lambda_{t_{N}}$ was defined in \eqref{Di}.

\begin{flushright}
$\square$
\end{flushright}

\subsection{Appendix B: Proof of Theorem 3.2}
We will prove the result for time consistent equilibrium; the proof for time inconsistent
equilibrium is similar and hence omitted. First consider the time period is $[(N-n)h,(N-n+1)h)$ and  choose an
arbitrary control $\pi=(\alpha,\beta)$ for any $J_{t_n}\in \{{\bf
0}, {\bf 1}  \}$ as follows:
\begin{eqnarray}
\pi=\left\{\begin{array}{ll} {\pi}^*_{t_n},
\ \ {\rm{for}}\ n=N-(n-1),N-(n-2)\cdots N-1,\\
\\\pi_{t_n},\ \ {\rm{for}}\ n=N-n.\label{alphan}\qquad\qquad\end{array} \right.
\end{eqnarray}
For convenience,  denote
$\varphi_{t_{N-n}}=\varphi(t_{N-n},C_{t_{N-n}},S_{t_{N-n}}).$ Assume $J_{t_{N-n}}=i.$ From
$$X_{t_N}^\pi=X^\pi_{t_{N-(n-1)}}+\sum_{k=N-(n-1)}^{N-1}\Delta
{X}_{t_k}^{*},$$ it follows that
\begin{eqnarray*}
&&\E_{t_{N-n}}[-\exp(-\gamma_{t_{N-n}}(X_{t_N}^\pi+I_{t_N})]\nonumber\\
&=&-e^{-\gamma_{t_{N-n}}x}\E_{t_{N-n}}\bigg[-e^{-\gamma_{t_{N-n}}\Delta
X_{t_{N-n}}^\pi}\cdot
\E_{t_{N-n+1}}[e^{-\gamma_{t_{N-n}} { [  {\displaystyle{\sum_{k=N-(n-1)}^{N-1}}}\Delta
{X}_{t_k}^{*}} + I_{t_N}]}]\bigg]\nonumber\\
 &=&
\E_{t_{N-n}}[-e^{-\gamma_{t_{N-n}}\Delta X^\pi_{t_{N-n}}}\cdot
e^{-\gamma_{t_{N-n}}Y^{*}_{t_{N-(n-1)}}}]\nonumber\\
&:=&-e^{-\gamma_{t_{N-n}}x}g_{N-n}(\alpha,\beta,\cdot).
\end{eqnarray*}
Here for  $k=N-(n-1),\cdots, N-2,N-1,$
$$\Delta
{X}_{t_k}^{*}=\alpha_{t_k}^*(\mu_{t_k}^ch+\sqrt{h}\sigma_{t_k}^c\Delta
b_{t_{k}}^1) +\Delta D_{t_k}+\varphi({t_k},C_{t_k},S_{t_k})h,$$ and
\begin{eqnarray}&&g_{N-n}(\alpha,\beta,\cdot)
= \E_{t_{N-n}}[-e^{-\gamma_{t_{N-n}}\Delta X^\pi_{t_{N-n}}}\cdot
e^{-\gamma_{t_{N-n}}Y^{*}_{t_{N-(n-1)}}}]\nonumber\\
&&=\frac{1}{2}e^{-\gamma_{t_{N-n}}\alpha(\mu^c_{t_{N-n}}h+\sigma^c_{t_{N-n}}\sqrt{h}
)}
 \E_{t_{N-n}}\bigg[e^{-\gamma_{t_{N-n}}\beta (\Delta D_{t_{N-n}}+\varphi_{t_{N-n}}h)}
 e^{-\gamma_{t_{N-n}}Y^{*}_{t_{N-n+1}}}| A_{t_{N-n}}\bigg]
\nonumber\\
&&+\frac{1}{2}e^{-\gamma_{t_{N-n}}\alpha(\mu^c_{t_{N-n}}h-\sigma^c_{t_{N-n}}\sqrt{h}
)}
 \E_{t_{N-n}}\bigg[e^{-\gamma_{t_{N-n}}\beta (\Delta D_{t_{N-n}}+\varphi_{t_{N-n}}h)}
 e^{-\gamma_{t_{N-n}}Y^{*}_{t_{N-n+1}}}| A_{t_{N-n}}^c\bigg]\nonumber\end{eqnarray} with
$A_{t_{N-n}}:=\{ \Delta b_{t_{N-n}}^1=1\}$ and
$A_{t_{N-n}}^c:=\{\Delta b_{t_{N-n}}^1=-1\}.$ Arguing as in the one
period case we get
$$g_{N-n}(0,0,\cdot)=\E_{t_{N-n}}[e^{-\gamma_{t_{N-n}}Y^{*}_{t_{N-n+1}}}]\leq\infty;$$
$$g_{N-n}(\infty,\beta,\cdot)=\infty;\ g_{N-n}(-\infty,\beta,\cdot)=\infty;$$
From arbitrage considerations it follows that
$$D_{t_{N-n}}(\omega)\in[\inf_{Q}\E^Q[D_{t_{N-n+1}}],
\sup_{Q}\E^Q[D_{t_{N-n+1}}]],$$ where $Q$ is the set of probability
measures. Thus
$$\Delta D_{t_{N-n}}(\omega)\in[D_{t_{N-n+1}}-\sup_{Q}
\E^Q[D_{t_{N-n+1}}],D_{t_{N-n+1}}-\inf_{Q}\E^Q[D_{t_{N-n+1}}]],$$ so the sets: $\{\omega:
\Delta D_{t_{N-n}}(\omega)>0\},$ and $\{\omega :\Delta
D_{t_{N-n}}(\omega)<0\}$ have positive probability. Consequently, it follows that:
$$g_{N-n}(\alpha,\infty,\cdot)=\infty;\ g_{N-n}(\alpha,-\infty,\cdot)=\infty.$$
Therefore the minimum of $g_{N-n}(\alpha,\beta,\cdot)$ is a
critical point. Hence
\begin{eqnarray}
\frac{\partial g_{N-n}}{\partial \alpha}
 &=&\frac{(\mu^c_{t_{N-n}}h+\sigma^c_{t_{N-n}}\sqrt{h}
)}{2}
e^{-\gamma_{t_{N-n}}\alpha(\mu^c_{t_{N-n}}h+\sigma^c_{t_{N-n}}\sqrt{h}
)}\nonumber\\&&\times
  \E_{t_{N-n}}\bigg[e^{-\gamma_{t_{N-n}}\beta^* (\Delta D_{t_{N-n}}+\varphi_{t_{N-n}}h)}
  e^{-\gamma(i)Y^{*}_{t_{N-n+1}}}
 | A_{t_{N-n}}\bigg]
\nonumber\\&+&\frac{(\mu^c_{t_{N-n}}h-\sigma^c_{t_{N-n}}\sqrt{h}
)}{2}
e^{-\gamma_{t_{N-n}}\alpha(\mu^c_{t_{N-n}}h-\sigma^c_{t_{N-n}}\sqrt{h}
)}\nonumber\\
&&\times
  \E_{t_{N-n}}\bigg[e^{-\gamma_{t_{N-n}}\beta^*(\Delta D_{t_{N-n}}+\varphi_{t_{N-n}}h)}
  e^{-\gamma_{t_{N-n}}Y^{*}_{t_{N-n+1}}}
| A_{t_{N-n}}^c\bigg]\nonumber\\
 &=&0.\end{eqnarray}

By direct calculation, we get that the optimal time consistent trading
strategy is
\begin{eqnarray}&&\alpha^*_{t_{N-n}}=\frac{1}{2\gamma_{t_{N-n}}\sigma_{t_{N-n}}^c\sqrt{h}}
\log[\frac{1+r_{t_{N-n}}^c\sqrt{h}}
{1-r_{t_{N-n}}^c\sqrt{h}}]\nonumber\\&+&
\frac{1}{2\gamma_{t_{N-n}}\sigma_{t_{N-n}}^c\sqrt{h}}\log(
\frac{\E_{t_{N-n}}[e^{-\gamma_{t_{N-n}}\beta_{t_{N-n}}^*
D_{t_{N-n+1}}}e^{-\gamma_{t_{N-n}}
Y^{*}_{t_{N-n+1}}}|A_{t_{N-n}}]}{\E_{t_{N-n}}[e^{-\gamma_{t_{N-n}}\beta_{t_{N-n}}^*
D_{t_{N-n+1}}}e^{-\gamma_{t_{N-n}}
Y^{*}_{t_{N-n+1}}}|A_{t_{N-n}}^c]}).\nonumber\\
\end{eqnarray} From the
equilibrium conditions it follows that $\beta_{t_{N-n}}^*=1.$ This combined with
$\frac{\partial g_{N-n}}{\partial \beta}=0,$ yield the equilibrium
price at $T_{N-n}.$ First, from $\frac{\partial g_{N-n}}{\partial
\beta}=0,$ one gets

\begin{eqnarray}
&& \E_{t_{N-n}} \bigg[
 (\Delta D_{t_{N-n}}+\varphi_{t_{N-n}}h)\cdot e^{-\gamma_{t_{N-n}}
 (\Delta D_{t_{N-n}}+\varphi_{t_{N-n}}h)}
 e^{-\gamma_{t_{N-n}}Y^{*}_{t_{N-n+1}}}| A_{t_{N-n}}\bigg]
\nonumber\\\!\!\!&=&\!\!\!-e^{2\gamma_{t_{N-n}}\alpha\sigma^c_{t_{N-n}}\sqrt{h}
}
 \E_{t_{N-n}} \bigg[
 (\Delta D_{t_{N-n}}+\varphi_{t_{N-n}}h)\cdot e^{-\gamma_{t_{N-n}}
 (\Delta D_{t_{N-n}}+\varphi_{t_{N-n}}h)}e^{-\gamma_{t_{N-n}}Y^{*}_{t_{N-n+1}}}|
 A_{t_{N-n}}^c\bigg]\nonumber.\end{eqnarray}
From $\frac{\partial g_{N-n}}{\partial \alpha}=0$ it follows that
$$e^{2\gamma_{t_{N-n}}\alpha^*\sigma^c_{t_{N-n}}\sqrt{h} }=
\frac{(1+r^c_{t_{N-n}}\sqrt{h})\E_{t_{N-n}} [e^{-\gamma_{t_{N-n}}\beta^*
D_{t_{N-n+1}}}e^{-\gamma_{t_{N-n}}Y^{*}_{t_{N-n+1}}}|{A_{t_{N-n}}
}]}{(1-r^c_{t_{N-n}}\sqrt{h})\E_{t_{N-n}} [e^{-\gamma(i)\beta^*
D_{t_{N-n+1}}}e^{-\gamma_{t_{N-n}}Y^{*}_{t_{N-n+1}}}|{ A_{t_{N-n}}^c}]}$$
This together with
$$\Delta D_{t_{N-n}}+\varphi_{t_{N-n}}h=D_{t_{N-n+1}}-\E^{Q^*}_{t_{N-n}}[D_{t_{N-n+1}}
],$$ 
 (here $Q^*$ is the equilibrium probability measure to be found) yield
\begin{eqnarray}
&&(\frac{2}{1-r^c_{t_{N-n}}\sqrt{h}})\E^{Q^*}_{t_{N-n}}[D_{t_{N-n+1}}
]
 \nonumber\\&=&
\frac{\E_{t_{N-n}}[D_{t_{N-n+1}} e^{-\gamma_{t_{N-n}}
D_{t_{N-n+1}}}e^{-\gamma_{t_{N-n}}Y^{*}_{t_{N-n+1}}}|A_{t_{N-n}}]}{\E_{t_{N-n}}[
e^{-\gamma_{t_{N-n}}
D_{t_{N-n+1}}}e^{-\gamma_{t_{N-n}}Y^{*}_{t_{N-n+1}}}|A_{t_{N-n}}]}\nonumber\\
&&+\frac{(1+r^c_{t_{N-n}}\sqrt{h})} {(1-r^c_{t_{N-n}}\sqrt{h})}
\frac{\E_{t_{N-n}}[D_{t_{N-n+1}} e^{-\gamma_{t_{N-n}}
D_{t_{N-n+1}}}e^{-\gamma_{t_{N-n}}
Y^{*}_{t_{N-n+1}}}|A_{t_{N-n}}^c]}{\E_{t_{N-n}}[ e^{-\gamma_{t_{N-n}}
D_{t_{N-n+1}}}e^{-\gamma_{t_{N-n}} Y^{*}_{t_{N-n+1}}}|A_{t_{N-n}}^c]}\nonumber.
\end{eqnarray}
Consequently

\begin{eqnarray}
&&D_{t_{N-n}}-\varphi_{t_{N-n}}h=   \frac{1-r_{t_n}^c\sqrt{h}}{2}  \frac{\E_{t_{N-n}}[D_{t_{N-n+1}} e^{-\gamma_{t_{N-n}}
D_{t_{N-n+1}}}e^{-\gamma_{t_{N-n}}Y^{*}_{t_{N-n+1}}}|A_{t_{N-n}}]}{\E_{t_{N-n}}[
e^{-\gamma_{t_{N-n}}
D_{t_{N-n+1}}}e^{-\gamma_{t_{N-n}}Y^{*}_{t_{N-n+1}}}|A_{t_{N-n}}]}\nonumber\\
&+&  \frac{1+r_{t_n}^c\sqrt{h}}{2}    \frac{\E_{t_{N-n}}[D_{t_{N-n+1}}
e^{-\gamma_{t_{N-n}}
D_{t_{N-n+1}}}e^{-\gamma_{t_{N-n}}Y^{*}_{t_{N-n+1}}}|A_{t_{N-n}}^c]}{\E_{t_{N-n}}[
e^{-\gamma_{t_{N-n}} D_{t_{N-n+1}}}e^{-\gamma_{t_{N-n}}Y^{*}_{t_{N-n+1}}}|A_{t_{N-n}}^c]}.\nonumber\\
\end{eqnarray}
Thus, the equilibrium price is
\begin{eqnarray}
D_{t_{N-n}}
&=&\E_{t_{N-n}}[(D_{t_{N-n+1}}+\varphi_{t_{N-n}}h)\Lambda^*_{t_{N-n+1}}],
\end{eqnarray}
with $\Lambda^*_{t_{N-n+1}}$ defined in \eqref{lambda}.

\subsection{Appendix C: Proof of Corollary 3.2}
We consider the time period $[(N-n)h,(N-n+1)h).$ Recall that
$$\hat{X}_{t_N}=x+\Delta
\hat{X}_{t_{N-n}}+{\displaystyle{\sum_{k=N-(n-1)}^{N-1}}}\Delta
\hat{X}_{t_k},$$
 where $$\Delta
\hat{X}_{t_k}=\hat{\alpha}_{t_k}(\mu_{t_k}^ch+\sqrt{h}\sigma_{t_k}^c\Delta
b_{t_{k}}^1) +\hat{\beta}_{t_k}( \Delta
D_{t_k}+\varphi(t_k,C_{t_k},S_{t_k})h),$$ for $ k=N-(n-1),\cdots,
N-2,N-1$. Thus, we have
\begin{eqnarray}&&\E_{t_{N-n}}(U'(\hat{W}_{t_N}))
\nonumber\\
&&= \E_{t_{N-n}}[-e^{-\gamma\Delta \hat{X}_{t_{N-n}}}\cdot
e^{-\gamma \hat{Y}_{t_{N-(n-1)}}}]\nonumber\\
&&=\frac{1}{2}e^{-\gamma\hat{\alpha}_{t_{N-n}}(\mu^c_{t_{N-n}}h+\sigma^c_{t_{N-n}}\sqrt{h}
)}e^{-\gamma\varphi_{t_{N-n}}h}
 \E_{t_{N-n}}\bigg[e^{-\gamma \Delta D_{t_{N-n}}}e^{-\gamma \hat{Y}_{t_{N-n+1}}}| A_{t_{N-n}}\bigg]
\nonumber\\
&&+\frac{1}{2}e^{-\gamma\hat{\alpha}_{t_{N-n}}(\mu^c_{t_{N-n}}h-\sigma^c_{t_{N-n}}\sqrt{h}
)}e^{-\gamma\varphi_{t_{N-n}}h}
 \E_{t_{N-n}}\bigg[e^{-\gamma \Delta D_{t_{N-n}}}
 e^{-\gamma \hat{Y}_{t_{N-n+1}}}| A_{t_{N-n}}^c\bigg].\nonumber\end{eqnarray} 

From direct calculations, it follows that on the set
$\{\omega:\omega\in A^c_{t_{N-n}}\}$
\begin{eqnarray}
&&\frac{\E_{t_{N-n+1}}[U'(\hat{W}_{t_N})]}{\E_{t_{N-n}}[{U'(\hat{W}_{t_N})}]}\nonumber\\
 &=&\frac{(1+r_{t_{N-n}}^c\sqrt{h})e^{\gamma\hat{\alpha}_{t_{N-n}}\sigma^c_{t_{N-n}}\sqrt{h}
}e^{-\gamma\hat{\alpha}_{t_{N-n}}\sigma^c\sqrt{h }\Delta b_{t_{N-n}^1}}e^{-\gamma
D_{t_{N-n+1}}}e^{-\gamma \hat{Y}_{t_{N-n+1}}}}{\E_{t_{N-n}}
\bigg[e^{-\gamma D_{t_{N-n+1}}}e^{-\gamma \hat{Y}_{t_{N-n+1}}}|
 A_{t_{N-n}}^c\bigg]}\nonumber\\
 &=&\lambda_{t_{N-n}}\frac{e^{-\gamma
D_{t_{N-n+1}}}e^{-\gamma \hat{Y}_{t_{N-n+1}}}}{\E_{t_{N-n}}[e^{-\gamma
D_{t_N-n+1}}e^{-\gamma \hat{Y}_{t_{N-n+1}}}|
 A_{t_{N-n}}^c]}\nonumber.
\end{eqnarray}
Moreover, on the set of $\{\omega:\omega\in A_{t_{N-n}}\}$
\begin{eqnarray}
\frac{\E_{t_{N-n+1}}[U'(\hat{W}_{t_N})]}{\E_{t_{N-n}}[{U'(\hat{W}_{t_N})}]}&=&
\lambda_{t_{N-n}}\frac{e^{-\gamma D_{t_{N-n+1}}}e^{-\gamma
\hat{Y}_{t_{N-n+1}}}}{\E_{t_{N-n}}[e^{-\gamma D_{t_N-n+1}}e^{-\gamma
\hat{Y}_{t_{N-n+1}}}|
 A_{t_{N-n}}]}\nonumber.
\end{eqnarray}
Therefore
\begin{eqnarray}
\frac{\E_{t_{N-n+1}}[U'(\hat{W}_{t_N})]}{\E_{N-n}[{U'(\hat{W}_{t_N})}]}=\hat{\Lambda}_{t_{N-n+1}}.
\end{eqnarray}


\begin{thebibliography}{99}
\bibitem{Bar}
{\sc Barberis, N. AND Huang, M.} (2001) {Mental Accounting, Loss
Aversion, and Individual Stock Returns}, {\em The Journal of
Finance} {\bf LVI 4}, 1247-1292.

\bibitem{Bhattacharya}
{\sc Bhattacharya, C., D.} (1981) {Variance Analysis to chances in
return of Investment}, Working Paper {\bf No. 365}, Indian Institute
of Management, Ahmedabad.



\bibitem{Bizid}
{\sc Bizid, A., AND Jouini, E.} {(2001)} {Incomplete markets and
short-sales constraints: an equilibrium approach}, {\em
International Journal of Theoretical and Applied Finance}, {\bf
4(2)},211-243.

\bibitem{Bjo}
{\sc Bj\"{o}rk, T. AND Murgoci, A.} (2010) {A General Theory of
Markovian Time Inconsistent Stochastic Control Problems}, {\em
Preprint}.

\bibitem{Brennan}
{\sc Brennan, M.} (1979) {The Pricing of Contingent Claims in
Discrete Time Models} {\em The Journal of Finance}, {\bf
34(1)},53-68.

\bibitem{Camara}
{\sc Camara, A.} (2003) {A Generalization of the Brennan-Rubinstein
Approach for the Pricing of Derivatives}, {\em The Journal of
Finance}, \textbf{VOL. LVIII, NO. 2},805-819.

\bibitem{Cao}
{\sc Cao, M., AND Wei, J.} (2004) {Pricing weather derivative: an
equilibrium approach}, {\em Journal of Futures Markets}, {\bf 24},
1065-1089.

\bibitem{Cherid}
{\sc Cheridito, P., Horst, U., Kupper, M., AND Pirvu, T., A.} (2011)
{Equilibrium Pricing in Incomplete Markets under Translation
Invariant Preferences},
http://papers.ssrn.com/sol3/papers.cfm?abstract


\bibitem{Danthine}
{\sc Danthine, J., Donaldsond, J., Giannikose, C., AND Guirguisf,
H.} (2004) {On the consequences of state dependent preferences for
the pricing of financial assets}, {\em Finance Research Letters}
{\bf 1(3)},143-153.

\bibitem{Davis}
{\sc Davis, M.} (2001) {Pricing weather derivatives by marginal
value}, {\em  Quantitative Finance}, {\bf 1}, 305-308.

\bibitem{EKE1}
{\sc Ekeland, I. AND Pirvu, T., A.} (2008) {Investment and
consumption without commitment}, {\em Mathematics and Financial
Economics,} {\bf 2}, 57-86.

\bibitem{Gordon}
{\sc Gordon, S., AND  St-Amour, P.} (2004) {Asset returns and
state-dependent risk preferences}, {\em Journal of Business and
Economic Statistics,} {\bf 22}, 241-252.


\bibitem{Gordon1}
{\sc Gordon, S., And  St-Amour, P.} (2000) {A preference regime model of bull and bear markets}, {\em American Economic Review,} {\bf 90}, 1019-1033.


\bibitem{Henderson0}
{\sc Henderson, V. } (2002) {Valuation of claims on non-traded assets using utility maximization,}
{\em Mathematical Finance,} {\bf 12}, 351-373. 




\bibitem{Henderson}
{\sc Henderson, V. AND Hobson, D.} (2004) {Utility indifference pricing-an overview,}
{\em Volume on Indifference Pricing}, (ed. R. Carmona), Princeton University Press. 



\bibitem{Hodges}
{\sc Hodges, S. AND Neuberger, A.} (1989) {Optimal replication of contingent claims under transaction costs,}
{\em Review future Markets}, {\bf  8}, 222-239. 


\bibitem{Horst}
{\sc Horst,U., Pirvu, T., A. AND Dos Reis, G.} (2010) {On
securitization, market completion and equilibrium risk transfer,}
{\em Mathematics and Financial Economics}, {\bf  2},211-252.

\bibitem{Lee1}
{\sc Lee, Y., AND Oren,S.}(2009) {An equilibrium pricing model for
weather derivatives in a multi-commodity setting}, {\em Energy
Economics,} {\bf 31}, 702-713.

\bibitem{Lee}
{\sc Lee, Y., AND Oren, S. } (2010) {A multi-period equilibrium
pricing model of weather derivatives}. {\em Energy System,} {\bf
1},3-30.



\bibitem{MZ0}
{\sc  Musiela, M. AND Zariphopoulou, T.} (2004)
{An example of indifference prices under exponential preferences}, {\em Finance and Stochastics} {\bf 8}, 229-239.




\bibitem{MZ1}
{\sc  Musiela, M. AND Zariphopoulou, T.} (2004)
{A valuation algorithm for indifference prices in incomplete markets}, {\em Finance and Stochastics} {\bf 8(3)}, 399-414.


\bibitem{MZ2}
{\sc  Musiela, M., Sokolova, E. AND Zariphopoulou, T.} (2009)
{Indiffeerence valuation in incomplete binomial 
models}.  \emph{Preprint} 





\bibitem{Rubi}
{\sc Rubinstein, M.} (1976) {The  Valuation  of Uncertain  Income
Streams and the Pricing of Options}, {\em Bell  Journal of
Economics,} {\bf 7}, 407-425.


\bibitem{Tra}
{\sc Pirvu, T.,A. AND Zhang, H.} (2011) {Utility Indifference
Pricing: A Time Consistent Approach,} {\em Submitted},
http://adsabs.harvard.edu/abs/2011arXiv1102.5075P

\bibitem{Luiz}
{\sc Vitiello, L. AND Poon. S., H.} (2010) {A General Equilibrium
and Preference Free Model for Pricing Options Under Transformed
Gamma Distribution}, {\em Journal of Futures Markets} {\bf 30(5)},
409-431.


\bibitem{Yuan}
{\sc Yuan, B. AND Chen, K.} (2006) {Impact of investors varying
risk aversion on the dynamics of asset price fluctuations}, {\em
Journal of Economic Interaction and Coordination,} {\bf 1(2)},
189-214.






\end{thebibliography}
\end{document}